\documentclass[reqno,12pt]{amsart}

\usepackage{epsf}
\usepackage{graphics}
\usepackage{amssymb}
\usepackage{amsmath}

\date{}

\theoremstyle{plain}
\newtheorem{theorem}{Theorem}
\newtheorem{corollary}{Corollary}

\theoremstyle{definition}

\theoremstyle{remark}

\newtheorem*{examples}{Examples}
\newtheorem*{remark}{Remark}

\title{Construction of Non-Alternating Knots} 

\author{Sebastian Baader}

\begin{document}

\begin{abstract} We investigate the behaviour of Rasmussen's invariant~$s$
under the sharp operation on knots and obtain a lower bound for the sharp
unknotting number. This bound leads us to an interesting move that transforms
arbitrary knots into non-alternating knots.
\end{abstract}

\maketitle

\section{Introduction}

An unknotting operation is a local operation that allows us to untie every
knot in finitely many steps. The most popular unknotting operation is a
simple crossing change. Every unknotting operation gives rise to a measure of
complexity for knots, called an unknotting number. An effective lower bound
for the usual unknotting number was introduced by Rasmussen (\cite{Ra}). His
invariant led to an easy computation of the genera and unknotting numbers of
torus knots. In this paper, we study the sharp unknotting operation via
Rasmussen's invariant~$s$. 

The sharp unknotting operation is a local move that acts on link diagrams, as
shown in figure~1. It has been introduced by Murakami (\cite{Mu}) and gives
rise to the unknotting number $u_\#$. The usual unknotting number is denoted
by $u$.

\begin{figure}[ht]
\scalebox{1}{\raisebox{-0pt}{$\vcenter{\hbox{\epsffile{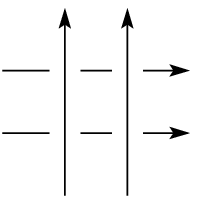}}}$}} \qquad
$\longleftrightarrow$ \qquad 
\scalebox{1}{\raisebox{-0pt}{$\vcenter{\hbox{\epsffile{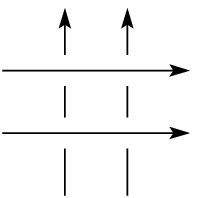}}}$}}
\caption{Sharp operation}
\end{figure}

Our main result involves a special sharp operation, called a positive sharp
operation. A positive sharp operation introduces eight positive crossings to a
link diagram, as shown in figure~2.

\begin{figure}[ht]
\scalebox{1}{\raisebox{-0pt}{$\vcenter{\hbox{\epsffile{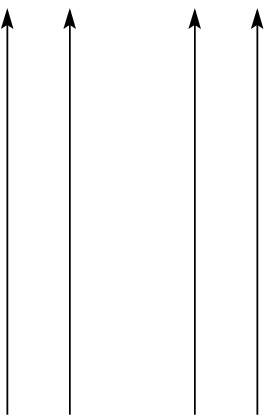}}}$}} \qquad
$\longrightarrow$ \qquad 
\scalebox{1}{\raisebox{-0pt}{$\vcenter{\hbox{\epsffile{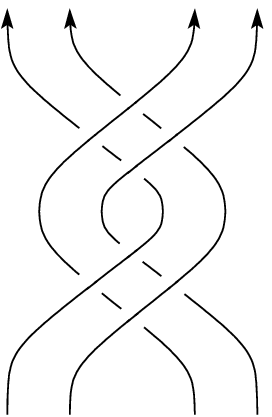}}}$}}
\caption{Positive sharp operation}
\end{figure}

For a diagram $D$ of
a knot $K$, we denote by $w(D)$, $O(D)$, and $g(K)$ the writhe of $D$ (i.e. the
algebraic crossing number of $D$), the number of Seifert circles of $D$, and
the genus of $K$, respectively.

\begin{theorem}
Let $D$ be any knot diagram, and suppose $D'$ is obtained from $D$ by the
application of~$n$ positive sharp operations. If 
$$n>g(K)-\frac{1+w(D)-O(D)}{2},$$
then $D'$ represents a non-alternating knot.
\end{theorem}

The quantity on the right hand side of the inequality in
theorem~1 is always positive, as follows from Bennequin's inequality
(\cite{Be}). For positive knot diagrams, it is actually zero.

\begin{corollary} Let $D$ be a positive knot diagram, and suppose $D'$ is
obtained from $D$ by the application of one positive sharp operation. Then $D'$
represents a non-alternating knot.
\end{corollary}

\begin{examples} \quad
\begin{enumerate}
\item [(i)] The closure of the braid $\sigma_1^{-1} \sigma_2 \sigma_1 \sigma_3
\sigma_2$ represents the trivial knot $O$. For the corresponding knot
diagram $D$ (see figure~3, on the left hand side), Bennequin's inequality is an equality:
$$g(O)-\frac{1}{2}(1+w(D)-O(D))=0-\frac{1}{2}(1+3-4)=0.$$
Therefore, if we apply one positive sharp operation at the top of this braid
diagram, we obtain a non-alternating knot. It is a 2-cable of the positive
trefoil knot.\\
\item [(ii)] The closure of the braid $\sigma_1 \sigma_2 \sigma_3$ represents
the trivial knot, too (see figure~3, on the right hand side). By corollary~1,
the application of one positive sharp operation at the top of that braid
diagram yields a non-alternating knot. This time, we obtain the knot
$10_{139}$, in Rolfsen's notation (\cite{Ro}).\\
\end{enumerate}
\end{examples}

\begin{figure}[ht]
\scalebox{1}{\raisebox{-0pt}{$\vcenter{\hbox{\epsffile{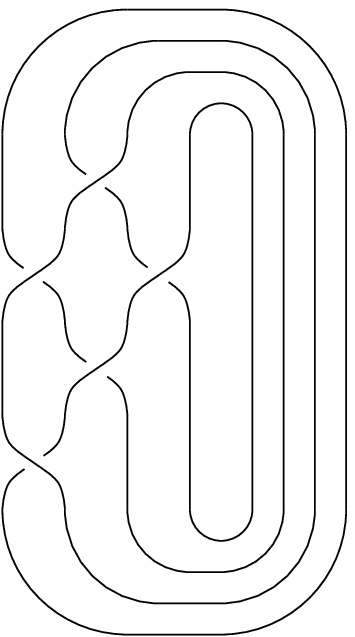}}}$}} \qquad
\qquad 
\scalebox{1}{\raisebox{-0pt}{$\vcenter{\hbox{\epsffile{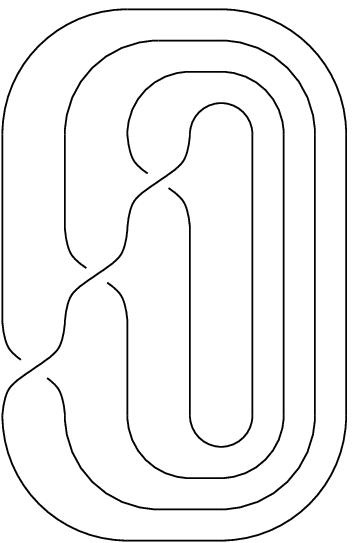}}}$}}
\caption{Two diagrams of the trivial knot}
\end{figure}

A sharp unknotting operation changes 4 crossings of a diagram. Therefore,
$u$ cannot exceed $4u_\#$. Rasmussen's invariant allows us to detect knots with $u(K)=4u_\#(K)$. 

\begin{theorem} \quad
\begin{enumerate}
\item $u_\#(K) \geqslant \frac{|s(K)|}{8}$.\\
\item If $u_\#(K)=\frac{|s(K)|}{8}$, then $K$ is either trivial or
non-alternating. In any case, $u(K)=4u_\#(K)$ holds.\\
\end{enumerate}
\end{theorem}

\begin{examples} (continued)
\begin{enumerate}
\item [(i)] The diagram $D$ of the 2-cable knot $K$ we constructed above has 4
Seifert circles and writhe 11, whence $1+w(D)-O(D)=8$. The latter quantity is
a lower bound for the invariant $s(K)$ 
(see \cite{Sh}). This proves $u_\#(K)=\frac{|s(K)|}{8}=1$. \\

\item [(ii)] The diagram $D$ of the knot $10_{139}$ we constructed above has 4
Seifert circles and writhe 11. Again, we conclude
$u_\#(10_{139})=\frac{|s(10_{139})|}{8}=1$.\\
\end{enumerate}
\end{examples}

\section{Rasmussen's Invariant and the Sharp Unknotting Operation}

The proofs of theorems~1 and~2 are based upon the following three properties
of Rasmussen's invariant~$s$:

\begin{enumerate}
\item $|s(K)| \leqslant 2u(K)$,\\
\item $s(K)=\sigma(K)$, for all alternating knots $K$ (here $\sigma(K)$ is the
signature of the knot $K$),\\
\item $1+w(D)-O(D) \leqslant s(K)$, where $D$ is any diagram of a knot $K$.\\
\end{enumerate}

The first two properties were proved by Rasmussen (\cite{Ra}), whereas the
third inequality was proved by Shumakovitch (\cite{Sh}). The main argument in
the proof of~(3) is Rudolph's reduction to the case of positive
diagrams (\cite{Ru}).

As we remarked after theorem~2, the usual unknotting number~$u$ cannot exceed
$4u_\#$. Together with the inequality~(1), this immediately proves the first
statement of theorem~2:
$$u_\#(K) \geqslant \frac{u(K)}{4} \geqslant \frac{|s(K)|}{8}.$$
In \cite{Mu}, Murakami proved the following estimate for $u_\#$, in terms of
the signature $\sigma$ of a knot:
$$u_\#(K) \geqslant \frac{|\sigma(K)|}{6}.$$
This implies the second statement of theorem~2:

Let $K$ be a knot with $u_\#(K)=\frac{|s(K)|}{8}$. Murakami's inequality tells
us that
$$\frac{|s(K)|}{8} \geqslant \frac{|\sigma(K)|}{6}.$$
If, in addition, $K$ is alternating, then $s(K)=\sigma(K)$, by~(2).
Therefore, $s(K)=\sigma(K)=0$, $u_\#(K)=0$, and $K$ is the trivial knot. In
any case, $4u_\#(K)=u(K)$ holds.

In order to prove theorem~1, we have to study the behaviour of the numbers
$w(D)$ and $O(D)$ under a positive sharp operation: a
positive sharp operation increases the writhe by 8 and leaves the number of
Seifert circles invariant. Now, let $D$ be any knot diagram of a knot $K$.
Further, suppose $D'$ is obtained from $D$ by the application of~$n$ positive
sharp operations. $D'$ represents a knot $K'$. Using~(3), we find the
following lower bound for $s(K')$:
$$s(K') \geqslant 1+w(D')-O(D')=1+w(D)+8n-O(D).$$
On the other hand, we have the following upper bound for the signature
$\sigma(K')$:
$$\sigma(K') \leqslant \sigma(K)+6n \leqslant 2g(K)+6n.$$
The first inequality is due to Murakami (\cite{Mu}): the signature of a knot
cannot increase by more than 6 under a sharp operation. The second inequality
is obvious, since the signature of a knot $K$ is the signature of a Seifert
matrix of size $2g(K)$. Now, if 
$$n>g(K)-\frac{1+w(D)-O(D)}{2},$$
then
$$s(K')-\sigma(K') \geqslant 2n-2g(K)+1+w(D)-O(D)>0,$$
whence $K'$ is non-alternating. This completes the proof of theorem~1.

\begin{remark} Throughout this paper, we could replace Rasmussen's invariant
$s$ by the concordance invariant $2\tau$ coming from knot Floer homology,
since the three properties (1), (2) and (3) are also valid for $2\tau$. A list
of properties that are shared by $s$ and $2\tau$ is contained in \cite{HO}. In
the same paper, M.~Hedden and Ph.~Ording show that the invariants $s$ and
$2\tau$ are not equal.
\end{remark}

\bigskip
\noindent
Department of Mathematics,
ETH Z\"urich, 
Switzerland

\bigskip
\noindent
\emph{sbaader@ethz.ch}

\end{document}